\documentclass[11pt,reqno,tbtags]{amsart}
\usepackage{amssymb}
\numberwithin{equation}{section}

\newtheorem{theorem}{Theorem}[section]
\newtheorem{lemma}[theorem]{Lemma}

\newtheorem{corollary}[theorem]{Corollary}

\theoremstyle{definition}

\newtheorem{remark}[theorem]{Remark}

\newtheorem*{ack}{Acknowledgement}

\theoremstyle{remark}

\newenvironment{romenumerate}{\begin{enumerate}% gives (i), (ii) etc.
 }{\end{enumerate}}

\newcommand{\refT}[1]{Theorem~\ref{#1}}
\newcommand{\refC}[1]{Corollary~\ref{#1}}
\newcommand{\refL}[1]{Lemma~\ref{#1}}
\newcommand{\refR}[1]{Remark~\ref{#1}}
\newcommand{\refS}[1]{Section~\ref{#1}}
\newcommand{\refand}[2]{\ref{#1} and~\ref{#2}}

\begingroup
  \divide\count255 by 60
  \multiply\count255 by -60
  \advance\count255 by \time
  \ifnum \count255 < 10 \xdef\klockan{\the\count1.0\the\count255}
  \else\xdef\klockan{\the\count1.\the\count255}\fi
\endgroup

\newcommand{\sumin}{\sum_{i=1}^n}

\newcommand\set[1]{\ensuremath{\{#1\}}}

\newcommand\xpar[1]{(#1)}
\newcommand\bigpar[1]{\bigl(#1\bigr)}
\newcommand\Bigpar[1]{\Bigl(#1\Bigr)}
\newcommand\biggpar[1]{\biggl(#1\biggr)}
\newcommand\lrpar[1]{\left(#1\right)}

\def\rompar(#1){\textup(#1\textup)}    % usage: \rompar(...)

\newcommand\parfrac[2]{\lrpar{\frac{#1}{#2}}}

\def\xexp(#1){e^{#1}}

\newcommand\ntoo{\ensuremath{{n\to\infty}}}
\newcommand\Ntoo{\ensuremath{{N\to\infty}}}

\newcommand\nutoo{\ensuremath{{\nu\to\infty}}}

\newcommand\ie{i.e.\spacefactor=1000}

\newcommand\cf{cf.\spacefactor=1000}

\newcommand{\tend}{\longrightarrow}
\newcommand\dto{\overset{\mathrm{d}}{\tend}}

\newcounter{CC} 
\newcounter{cc}
\newcommand\E{\operatorname{\mathbb E{}}}
\renewcommand\P{\operatorname{\mathbb P{}}}

\newcommand\Po{\operatorname{Po}}

\newcommand\fall[1]{^{\underline{#1}}}

\newcommand\ga{\alpha}
\newcommand\gb{\beta}
\newcommand\gd{\delta}

\newcommand\gl{\lambda}
\newcommand\gL{\Lambda}

\newcommand\cA{\mathcal A}
\newcommand\cB{\mathcal B}

\newcommand\cE{\mathcal E}
\newcommand\cH{\mathcal H}

\newcommand\tW{{\widetilde W}}
\newcommand\tY{{\widetilde Y}}

\newcommand\ett[1]{\boldsymbol1[#1]} 

\def\[#1]{[\![#1]\!]}

\newcommand\qq{^{1/2}}
\newcommand\qqw{^{-1/2}}
\newcommand\qw{^{-1}}

\renewcommand{\=}{:=}

\newcommand\dtv{d_{\mathrm{TV}}}

\newcommand{\pgf}{probability generating function}

\newcommand\lhs{left-hand side}
\newcommand\rhs{right-hand side}

\newcommand\gnd{\ensuremath{G(n,(d_i)_1^n)}}
\newcommand\gndx{\ensuremath{G^*(n,(d_i)_1^n)}}
\newcommand\nn{^{(n)}}

\newcommand\dn{\ensuremath{(d_i)_1^n}}
\newcommand\dnn{\ensuremath{(d_i\nn)_1^n}}
\newcommand\hdn{\ensuremath{(\hd_i)_1^{\hn}}}
\newcommand\hd{\hat{d}}
\newcommand\hn{\hat{n}}
\newcommand\nnu{^{(\nu)}}
\newcommand\sumi{\sum_i}
\newcommand\sumd{\sum_i d_i}
\newcommand\sumdd{\sum_i d_i^2}
\newcommand\sumhdd{\sum_i \hd_i^2}
\newcommand\sumv{\sum_v}

\newcommand\dbound{$\sum_v d_v^2=O(N)$}
\newcommand\maxd{\max_i d_i}
\newcommand\gnu{G^*_\nu}
\newcommand\gx{G^*}
\newcommand\gxk{G^*_k}
\newcommand\gxki{G^*_{k-1}}
\newcommand\tgxk{{\overline G}^*_k}
\newcommand\hgx{\widehat G^*}
\newcommand\simple{\text{ is simple}}
\newcommand\ii{_{ii}}
\newcommand\ij{_{ij}}
\newcommand\vw{_{vw}}
\newcommand{\prodil}{\prod_{i=1}^l}
\newcommand{\prodjm}{\prod_{j=1}^m}
\newcommand{\prodilx}[1]{\prod_{i=#1}^l}
\newcommand{\prodjmx}[1]{\prod_{j=#1}^m}
\newcommand\hu[1]{u^{(#1)}}
\newcommand\hv[1]{v^{(#1)}}
\newcommand\hw[1]{w^{(#1)}}
\newcommand\bI{\bar I}
\newcommand\bIu{\bI_u}
\newcommand\bJ{\bar J}
\newcommand\bJe{\bJ_e}
\newcommand\bY{\overline Y}
\newcommand\dvh{d_{v;H}}
\newcommand\hh{G}
\newcommand\CS{Cauchy--Schwarz}
\newcommand\CSineq{\CS{} inequality}

\newcommand\REM[1]{{\raggedright\texttt{[#1]}\par\marginal{XXX}}}

\newcommand\urladdrx[1]{{\urladdr{\def~{{\tiny$\sim$}}#1}}}
\begin{document}
\title%[]
{The probability that a random multigraph is simple}

\date{September 28, 2006}

\author{Svante Janson}
\address{Department of Mathematics, Uppsala University, PO Box 480,
SE-751~06 Uppsala, Sweden}
\email{svante.janson@math.uu.se}
\urladdrx{http://www.math.uu.se/~svante/}
%\urladdr{http://www.math.uu.se/\~{}svante/}

%\keywords{}
%%{Primary: <subject>; Secondary: <subject>}

\begin{abstract} 
Consider a random multigraph $G^*$ with given vertex degrees
$d_1,\dots,d_n$,
contructed by the configuration model.
We show that, asymptotically for a sequence of such multigraphs with
the number of edges $\tfrac12\sumd\to\infty$,
the probability that the multigraph is simple stays away from 0 if and
only if $\sumdd=O\bigpar{\sumd}$. This was previously known only under
extra assumtions on the maximum degree $\maxd$.
We also give an asymptotic formula for this probability, extending
previous results by several authors.
\end{abstract}

\maketitle

\section{Introduction}\label{S:intro}

If $n\ge1$ and $\dn$ is a sequence of non-negative integers, we let
\gnd{} be the random (simple) graph with the $n$ vertices $1,\dots,n$,
and with vertex degrees $d_1,\dots,d_n$, uniformly chosen among all
such graphs (provided that there are any such graphs at all; in particular,
$\sum_i d_i$ has to be even). 
A standard method to study $\gnd$ is to consider the related random
multigraph \gndx{} defined by taking a set of $d_i$ \emph{half-edges}
at each vertex $i$ and then joining the 
half-edges into edges by taking a random partition of the set of all
half-edges into pairs; see \refS{Sprel} for details. 
This is known as the configuration model, 
and such a partition of the half-edges is known as a \emph{configuration};
this was
introduced by 
Bollob\'as \cite{BB80},
see also Section II.4 of \cite{bollobas}.
(See Bender and Canfield
\cite{BenderC} and Wormald \cite{WormaldPhD,Wormald81} for related
arguments.)

Note that \gndx{} is defined for all $n\ge1$ and all sequences $\dn$
such that $\sumd$ is even (we tacitly assume this throughout the
paper), 
and that we obtain \gnd{} if we condition \gndx{} on being a simple graph. 
The idea of using the configuration method to study \gnd{}
is that \gndx{} in many respects
is a simpler object than \gnd; 
thus it is often possible to show
results for \gnd{} by first  
studying \gndx{} and then conditioning on this multigraph being
simple. 
It is then of crucial importance to be able to estimate the probability
that \gndx{} is simple, and in particular to decide whether 
\begin{equation}
\label{a2}
\liminf_{\ntoo}\P\bigpar{\gndx\text{ is simple}}>0  
\end{equation}
for  given sequences $\dn=(d_i\nn)_1^n$ (depending on $n\ge1$).
(Note that \eqref{a2} implies that any statement holding for \gndx{}
with probability tending to 1 does so for \gnd{} too.) 

A natural condition that has been used by several authors 
using the configuration method (including myself \cite{SJ184})
as a sufficient condition for \eqref{a2} is
\begin{equation}
\label{a3}
\sumin d_i=\Theta(n)
\qquad\text{and}\qquad
\sumin d_i^2=O(n)
\end{equation}
together with some bound on $\maxd$.
(Recall that $A=\Theta(B)$ means that both $A=O(B)$ and $B=O(A)$
hold.) Results showing, or implying, 
that \eqref{a3} and a condition on $\maxd$ imply \eqref{a2}
have also been given by
several authors, for example
Bender and Canfield \cite{BenderC} with $\maxd=O(1)$;
Bollob\'as \cite{BB80}, see also Section II.4 in \cite{bollobas},
with
$\maxd\le\sqrt{2\log n}-1$; 
McKay \cite{McKay} with $\maxd=o(n^{1/4})$;
McKay and Wormald \cite{McKayWo} with $\maxd=o(n^{1/3})$.
(Similar results have also been proved for
bipartite graphs \cite{McKayrect},
digraphs \cite{cooperF-digraph},
and hypergraphs \cite{cooper-hyper}.)

Indeed, it is not difficult to see that 
the method used by Bollob\'as \cite{BB80,bollobas}  works, assuming
\eqref{a3}, provided only $\maxd=o(n\qq)$, see \refS{Spoisson}. This
has undoubtedly been noted by several experts, but we have not been
able to find a reference to it in print when we have needed one.

One of our main result is that, in fact, \eqref{a3} is sufficient for
\eqref{a2} without any assumption on $\maxd$, even in cases where the
Poisson approximation fails. Moreover, \eqref{a3} is essentially
necessary.

We remark that 
several papers (including several of the references given above) 
study $\P(\gndx\simple)$ from another point of view,
namely by studying the number of simple graphs with given degree sequence
$\dn$. It is easy to count configurations, and it follows that this
number equals, with $N$ the number of edges, see \eqref{N} below,
\begin{equation*}
  \frac{(2N)!}{2^N N!\, \prod_i d_i!} \P\bigpar{\gndx\simple};
\end{equation*}
such results are thus equivalent to results for
$\P\bigpar{\gndx\text{ is simple}}$. However, in this setting it is
also interesting to obtain detailed asymptotics
when $\P\bigpar{\gndx\simple}\to0$;
such results are included in several of the references above, 
but will not be treated here.

We will throughout the paper let $N$ be the number of edges in
\gndx. Thus
\begin{equation}\label{N}
  2N=\sumin d_i .
\end{equation}
It turns out that it is more natural to state our results in terms of
$N$ than $n$ (the number of vertices). We can state our first result
as follows; we use an index $\nu$ to emphasize that the result is
asymptotic, and thus should be stated for a sequence (or another
family) of multigraphs.

\begin{theorem}
  \label{T1}
Consider a sequence of random multigraphs 
$\gnu=G^*\bigpar{n_\nu,(d_i\nnu)_1^n}$. Let $N_\nu=\tfrac12\sumd\nnu$,
the number of edges in $\gnu$, and assume that, as $\nu\to\infty$,
$N_\nu\to\infty$. Then
\begin{romenumerate}
  \item
$\liminf_\nutoo \P(\gnu\simple)>0$ if and only if
$\sum_i (d_i\nnu)^2 = O(N_\nu)$;
  \item
$\lim_\nutoo\P(\gnu\simple)=0$ if and only if
$\sum_i (d_i\nnu)^2/N_\nu\to\infty$.
\end{romenumerate}
\end{theorem}

In the sequel we will for simplicity omit the index $\nu$, but all
results should be interpreted in the same way as \refT{T1}.

Usually, one studies $\gndx$ as indexed by $n$. We then have the
following special case of \refT{T1}, which includes the claim above
that \eqref{a3} is sufficient for \eqref{a2}.

\begin{corollary}
  \label{C1}
  Let $\dn=\dnn$ be given for $n\ge1$. Assume that
  $N=\Theta(n)$. Then, as \ntoo,
\begin{romenumerate}
  \item
$\liminf_\ntoo \P\bigpar{\gndx\simple}>0$ if and only if
$\sum_i (d_i\nnu)^2 = O(n)$,
  \item
$\P\bigpar{\gndx\simple}\to0$ if and only if
$\sum_i (d_i\nnu)^2/n\to\infty$.
\end{romenumerate}
\end{corollary}

\begin{remark}
  \label{R1}
Although we have stated \refC{C1} as a special case of \refT{T1} with
$N=O(n)$, it is essentially equivalent to \refT{T1}.
In fact, we may ignore all vertices of degree 0; thus we may assume
that $d_i\ge1$ for all $i$, and hence $2N\ge n$. If further $\sum_i
d_i^2=O(N)$, the \CSineq{} yields
\begin{equation*}
  2N=\sumin d_i
\le \Bigpar{n\sumin d_i^2}\qq = O(\sqrt{nN}),
\end{equation*}
and thus $N=\Theta(n)$.
In the case $\sumdd/n\to\infty$, it is possible to reduce some $d_i$
to 1 such that then $N=\tfrac12\sumd=\Theta(n)$ and still
$\sumdd/N\to\infty$; we omit the details since our proof does not use
this route.
\end{remark}

Our second main result is an asymptotic formula for the probability
that \gndx{} is simple.
\begin{theorem}
  \label{T2}
Consider \gndx{} and assume that $N\=\tfrac12\sumd \to\infty$.
Let $\gl\ij\=\sqrt{d_i(d_i-1)d_j(d_j-1)}/(2N)$; in particular
$\gl\ii=d_i(d_i-1)/(2N)$. Then 
\begin{multline}\label{t2a}
\P\bigpar{\gndx\simple}
\\
=
\exp\Bigpar{-\tfrac12\sum_i\gl\ii
 -\sum_{i<j}\bigpar{\gl\ij-\log(1+\gl\ij)}}
+o(1);
\end{multline}
equivalently,
\begin{multline}\label{t2b}
\P\bigpar{\gndx\simple}
\\
=
\exp\biggl(-\frac14\Bigpar{\frac{\sumi d_i^2}{2N}}^2+\frac14
+ \frac{\sumi d_i^2(d_i-1)^2}{16N^2}
+\sum_{i<j}\bigpar{\log(1+\gl\ij)-\gl\ij+\tfrac12\gl\ij^2}\biggr)
\\+o(1).
\end{multline}
\end{theorem}

In many cases, $\sumi d_i^2(d_i-1)^2$ in \eqref{t2b} may be replaced by
the simpler $\sumi d_i^4$; for example, this can be done whenever
\eqref{a2} holds, by \refT{T1} and \eqref{sofie}. Note, however, that
this is not always possible; a trivial counter example is obtained
with $n=1$ and $d_1=2N\to\infty$.

In the case $\maxd=o(N\qq)$, \refT{T2} simplifies as follows; see
also \refS{Spoisson}.
\begin{corollary}
  \label{C2}
Assume that $\Ntoo$ and $\maxd=o(N\qq)$.
Let
\begin{equation}\label{Lambda}
 \gL\=\frac1{2N}\sumin\binom{d_i}2=\frac{\sumi d_i^2}{4N}-\frac12.
\end{equation}
Then
\begin{equation*}
  \begin{split}
\P\bigpar{\gndx\simple}
&=
\exp\bigpar{-\gL-\gL^2}+o(1)
\\
&=
\exp\Bigpar{-\frac14\Bigpar{\frac{\sumi d_i^2}{2N}}^2+\frac14}
+o(1).
  \end{split}
\end{equation*}
\end{corollary}
This formula is well known, at least under stronger conditions on
$\maxd$, see, for example, 
Bender and Canfield \cite{BenderC}, 
Bollob\'as \cite[Theorem II.16]{bollobas}, 
McKay \cite{McKay} and McKay and Wormald \cite[Lemma 5.1]{McKayWo}.

\section{Preliminaries}\label{Sprel}

We introduce some more notation. 

We will often write $\gx$ for the random multigraph $\gndx$.

Let $V_n=\set{1,\dots,n}$; this is the vertex set of $\gndx$.
We will in the sequel denote elements of $V_n$ by $u,v,w$, possibly
with indices.
$V_n$ is also the vertex set of the complete graph $K_n$, and we let
$E_n$ denote the edge set of $K_n$; thus $E_n$ consists of the 
$\binom n2$ unordered pairs \set{v,w}, with $v,w\in V_n$ and $v\neq w$.
We will use the notation $vw$ for the edge $\set{v,w}\in E_n$.

For any multigraph $G$ with vertex set $V_n$, and $u\in V_n$, 
we let
$X_u(G)$ be the number of loops at $u$.
Similarly, if $e=vw\in E_n$, we let $X_e(G)=X\vw(G)$ be the number of
edges between $v$ and $w$.
We define further the indicators
\begin{equation*}
\begin{aligned}
  I_u(G)&\=\ett{X_u(G)\ge1},\quad && u\in V_n,
\\
  J_e(G)&\=\ett{X_e(G)\ge2}, && e\in E_n,
\end{aligned}  
\end{equation*}
and their sum
\begin{equation}\label{y}
  Y(G)\=\sum_{u\in V_n}I_u(G)+\sum_{e\in E_n}J_e(G).
\end{equation}
Thus $G$ is a simple graph if and only if $Y(G)=0$, and our task is to
estimate $\P(Y(\gx)=0)$.

As said above, the idea of the configuration model is that
we fix a set of $d_v$ half-edges for every vertex $v$; 
we denote these half-edges by $\hv{1},\dots,\hv{d_v}$, 
and say that they belong to $v$, or are at $v$. 
These sets are assumed to be disjoint, so the total number of
half-edges is $\sumv d_v=2N$. A configuration is a partition of the
$2N$ half-edges
into $N$ pairs, and each configuration defines a multigraph with
vertex set $V_n$ and vertex degrees $d_v$ by letting 
every pair $\set{x,y}$ of half-edges
in the configuration define an edge; if $x$ is a half-edge at $v$
and $y$ is a half-edge at $w$, we form an edge between $v$ and $w$
(and thus a loop if $v=w$). We express this construction by saying
that we join the two half-edges $x$ and $y$ to an edge; we may denote
this edge by $xy$.
Recall that $\gx$ is the random multigraph obtain 
from a (uniform) random configuration
by this construction.

We will until \refS{Sproofs} assume that
\begin{equation}\label{dbound}
  \sumv d_v^2 = O(N),
\end{equation}
\ie, that $\sumv d_v^2\le CN$ for some constant $C$. (The constants
implicit in the estimates below may depend on this constant $C$.)
Note that an immediate consequence is
\begin{equation}\label{dmaxbound}
  \max_v d_v = O(N\qq)=o(N).
\end{equation}
We may thus assume that $N$ is so large that $\max_v d_v< N/10$, say,
and thus all terms like $N-d_v$ are of order $N$. (The estimates we
will prove are trivially true for any finite number of $N$ by taking
the implicit constants large enough; thus it suffices to prove them
for large $N$.)

Note further that \eqref{dbound} implies, using \eqref{dmaxbound},
that for any fixed $k\ge 2$ 
\begin{equation}\label{sofie}
  \sumv d_v^k 
\le(\max_v d_v)^{k-2} \sumv d_v^2
= O(N^{k/2}).
\end{equation}

We further note that we can assume $d_v\ge1$ for all $v$, since
vertices with degree $0$ may be removed without any difference to our
results. (This is really not necessary, but it means that we do not
even have to think about, for example, $d_v\qw$ in some formulas below.)

We will repeatedly use the \emph{subsubsequence principle}, which says
that if $(x_n)_n$ is a sequence of real numbers and $a$ is a number such
that every subsequence of $(x_n)_n$ has a subsequence that converges to
$a$, then the full sequence converges to $a$. (This holds in any
topological space.) 

We denote the falling factorials by $x\fall k\=x(x-1)\dotsm(x-k+1)$.

\section{Two probabilistic lemmas}\label{SP}

We will use two simple probabilistic lemmas. The first is (at least
part (i)) a standard extension of the inclusion-exclusion principle;
we include a proof for completeness.

\begin{lemma}
  \label{LP}
Let $W$ be a non-negative integer-valued random variable such that 
$\E R^W<\infty$ for some $R>2$. 
  \begin{romenumerate}
\item
Then, 
for every $j\ge0$,
\begin{equation*}
\P(W=j) 
= \sum_{k=j}^\infty (-1)^{k-j}\binom kj \frac{1}{k!} \E(W\fall k).
\end{equation*}
\item
More generally, for every random variable $Z$ such that
$\E(|Z|R^W)<\infty$ for some $R>2$, and every $j\ge0$,
\begin{equation*}
  \E(Z\cdot\ett{W=j}) 
= \sum_{k=j}^\infty (-1)^{k-j}\binom kj \frac{1}{k!} \E(Z W\fall k).
\end{equation*}
  \end{romenumerate}
\end{lemma}

\begin{proof}
For (i),  let $f(t)\=\E (t^W)=\sum_j\P(W=j)t^j$ 
be the \pgf{} of $W$; this is by assumption convergent for $|t|\le R$,
at least. If $|t|\le R-1$ we have
\begin{equation*}
  f(t+1)=\E(1+t)^W=\sum_{k=0}^\infty \E \binom Wk t^k
\end{equation*}
and thus, if $|t|\le R-2$,
\begin{equation}\label{manne}
  f(t)
=\sum_{k=0}^\infty \E \binom Wk (t-1)^k
=\sum_{k=0}^\infty \E\bigpar{W\fall k/k!}
 \sum_{j=0}^\infty \binom kj t^j(-1)^{j-k}.
\end{equation}
The double series is absolutely convergent since
\begin{equation*}
\sum_{k=0}^\infty  \sum_{j=0}^\infty
\E\bigpar{ W\fall k/k!} \binom kj |t|^j
=\sum_{k=0}^\infty \E \binom Wk (|t|+1)^k
=f(|t|+2)<\infty.
\end{equation*}
Hence the result follows by extracting the coefficients of $t^j$ in
\eqref{manne} 

Part (ii) is proved it the same way, using instead $f(t)\=\E(Zt^W)$.
\end{proof}

The next lemma could be proved by \refL{LP} if made the hypothesis
somewhat stronger, but we prefer another proof.

\begin{lemma}
  \label{LB}
Let $(W_\nu)_\nu$ and $(\tW_\nu)_\nu$ be two sequences of non-negative
integer-valued random variables such that, for some $R>1$
\begin{equation}
  \label{lb1}
\sup_\nu \E \bigpar{R^{W_\nu}} <\infty
\end{equation}
and, for each fixed $k\ge1$,
\begin{equation}
  \label{lb2}
\E\bigpar{W_\nu\fall k} - \E\bigpar{\tW_\nu\fall k} \to 0
\qquad \text{as \nutoo}.
\end{equation}
Then, as \nutoo,
\begin{equation}
  \label{lb3}
\P(W_\nu=0)-\P(\tW_\nu=0)\to0.
\end{equation}
\end{lemma}

\begin{proof}
By the subsubsequence principle,
  it suffices to prove that every subsequence has a subsequence along
  which \eqref{lb3} holds. Since \eqref{lb1} implies that the sequence
$(W_\nu)_\nu$  is tight, we can by selecting a suitable subsequence
  assume that $W_\nu\dto W$ for some random variable $W$
(see Sections 5.8.2 and 5.8.3 in Gut \cite{Gut}). 
Moreover,
  \eqref{lb1} implies uniform integrability of the powers $W_\nu^k$
  for each $k$, and we thus have, as \nutoo{} along the selected subsequence,
$\E (W_\nu^k)\to \E (W^k)$ for every $k$ and thus also 
$\E \bigpar{W_\nu\fall k}\to\E \bigpar{W\fall k}$
(see Theorems 5.4.2 and 5.5.9 in \cite{Gut}). 
By \eqref{lb2}, this yields also
$\E \bigpar{\tW_\nu\fall k}\to\E \bigpar{W\fall k}$.
Furthermore, \eqref{lb1} implies by Fatou's lemma
(Theorem 5.5.8 in \cite{Gut}) that
$\E (R^W) \le \liminf \E (R^{W_\nu})<\infty$, or 
$\E\bigpar{e^{tW}}<\infty$ with $t=\log R>0$; hence the distribution of $W$ is
  determined by its moments (see Section 4.10 in \cite{Gut}).
Consequently, by the method of moments (Theorem 6.7 in \cite{JLR}),
still along the subsequence,
$\tW_\nu\dto W$ and thus
  \begin{equation*}
\P(W_\nu=0)-\P(\tW_\nu=0)\to \P(W=0)-\P(W=0)=0.	
  \end{equation*}
\vskip-\baselineskip
\end{proof}

\begin{remark}\label{RB}
  The same proof gives the stronger statement
$$
\dtv(W_\nu,\tW_\nu)\=\sum_j|\P(W_\nu=j)-\P(\tW_\nu=j)|\to0.
$$
\end{remark}

\section{Individual probabilities}\label{SB}

We begin by estimating the probabilities $\P(I_u(\gx)=1)$
and $\P(J\vw(\gx)=1)$. The following form will be convenient.

\begin{lemma}
  \label{L1}
Suppose \dbound. Then, for $\gx$, and for all $u,v,w\in V_n$,
if $N$ is so large that $d_u\le N$,
\begin{align*}
-\log\P(I_u=0)
&=-\log\P(X_u=0)=\frac{d_u(d_u-1)}{4N}+O\parfrac{d_u^3}{N^2}
\intertext{and, with 
  $\gl\vw\=\sqrt{d_v(d_v-1)d_w(d_w-1)}/(2N)$ as in \refT{T2},}
-\log\P(J\vw=0)
&=-\log\P(X\vw\le1)
\\
&=
-\log(1+\gl\vw)+\gl\vw+ O\parfrac{(d_v+d_w)d_v^2d_w^2}{N^3}.
\end{align*}
\end{lemma}

\begin{proof}
The calculation for loops is simple. We construct the random
configuration by first choosing partners to the half-edges at $u$,
one by one. A simple counting shows that
\begin{equation*}
\P(X_u=0)=\prod_{i=1}^{d_u}	\lrpar{1-\frac{d_u-i}{2N-2i+1}}
\end{equation*}
and thus, for large $N$, using $-\log(1-x)=x+O(x^2)$ when $|x|\le1/2$,
\begin{equation*}
-\ln\P(X_u=0)=\sum_{i=1}^{d_u}
 \lrpar{\frac{d_u-i}{2N}+O\parfrac{d_u^2}{N^2}}
=\frac{d_u(d_u-1)}{4N}+O\parfrac{d_u^3}{N^2}.
\end{equation*}

For multiple edges, a similar direct approach would be much more
complicated because of the possibility of loops at $v$ or $w$.
We instead use \refL{LP}(i), with $W=X\vw$.
We may assume $d_v,d_w\ge2$, since the result otherwise is trivial.
$X\vw\fall k$ is the number of ordered $k$-tuples of edges between $v$
and $w$; the corresponding pairs in the configuration may be chosen in
$d_v\fall kd_w\fall k$ ways, and each such set of $k$ pairs appears in
the configuration with probability
$((2N-1)(2N-3)\dotsm(2N-2k+1))\qw$. Thus
\begin{equation}\label{magnus}
  \E X\vw\fall k 
=
\frac{d_v\fall kd_w\fall k}{(2N-1)\dotsm(2N-2k+1)}
=
\frac{d_v\fall kd_w\fall k}{2^k(N-1/2)\fall k}
=
\prod_{i=0}^{k-1}\frac{(d_v-i)(d_w-i)}{2(N-1/2-i)}.
\end{equation}
In particular, 
\begin{equation}\label{jesper}
\E X\vw\fall2=\gl\vw^2(1+O(1/N))  
\end{equation}
and if $k\ge 3$, uniformly in $k$,
\begin{equation}\label{b2}
  \E X\vw\fall k 
=
(1+O(1/N))  \gl\vw^2
\prod_{i=2}^{k-1}\frac{(d_v-i)(d_w-i)}{2(N-1/2-i)}.
\end{equation}

Since $d_v<N$ (for large $N$ at least, see \eqref{dmaxbound}), the
ratios $(d_v-i)/(N-1/2-i)$ decrease as $i$ increases; hence, for large
$N$ and $i\ge2$,
\begin{equation*}
  \frac{d_v-i}{2(N-i-1/2)}
\le
  \frac{d_v-2}{2N-5}
<
  \frac{d_v-1}{2N}
<
  \frac{\sqrt{d_v(d_v-1)}}{2N}
\end{equation*}
and \eqref{b2} yields, uniformly for $k\ge3$,
\begin{equation}\label{b2b}
  \E X\vw\fall k 
\le
(1+O(1/N))  \gl\vw^2 \prod_{i=2}^{k-1}\gl\vw
= 
(1+O(1/N))  \gl\vw^k.
\end{equation}
In the opposite direction, still uniformly for $k\ge3$,
\begin{multline*}
  %\begin{equation*}
  %\begin{split}
  \E X\vw\fall k /\gl\vw^k
\ge
\prod_{i=2}^{k-1}\frac{(d_v-i)(d_w-i)}{2N\gl\vw}
\ge
\prod_{i=2}^{k-1}\frac{(d_v-i)(d_w-i)}{d_vd_w}
\\
=
\prod_{i=2}^{k-1}(1-i/d_v)(1-i/d_w)
\ge
1-\sum_{i=2}^{k-1}(i/d_v+i/d_w)
\ge 1-\frac{k^2}{d_v}-\frac{k^2}{d_w}.	
%  \end{split}
%\end{equation*}
\end{multline*}
Together with \eqref{b2b}, this shows that, uniformly for $k\ge 3$,
\begin{equation}\label{emma}
  \E X\vw\fall k 
=\gl\vw^k\lrpar{1+O\lrpar{k^2(d_v\qw+d_w\qw)}}.	
\end{equation}

We now use \refL{LP}(i), noting that trivially 
$\E R^{X\vw}\le R^{d_v}<\infty$ for every $R$. 
Thus, using \eqref{jesper} and
\eqref{emma} and observing that  $\gl\vw=O(1)$ by
\eqref{dmaxbound},
\begin{equation*}
  \begin{split}
\P(&X\vw=0)
=1-\E X\vw + \tfrac12\gl\vw^2(1+O(N\qw))
\\
&\hskip 12em
+\sum_{k=3}^\infty(-1)^k\frac{\gl\vw^k	}{k!}
\lrpar{1+O\lrpar{k^2(d_v\qw+d_w\qw)}}
\\&
=
1-\E X\vw + e^{-\gl\vw}-1+\gl\vw
  +O\lrpar{\gl\vw^2N\qw+\gl\vw^3e^{\gl\vw}(d_v\qw+d_w\qw)}
\\&
=
e^{-\gl\vw}+\gl\vw-\E X\vw +O\lrpar{(d_v^2d_w^3+d_v^3d_w^2)N^{-3}}
.
  \end{split}
\end{equation*}
Similarly, still by \refL{LP}(i), 
\begin{equation*}
  \begin{split}
\P(X\vw&=1)
=\E X\vw - \gl\vw^2(1+O(N\qw))
\\
&\hskip 8em
+\sum_{k=3}^\infty(-1)^{k-1}\frac{\gl\vw^k	}{(k-1)!}
\lrpar{1+O\lrpar{k^2(d_v\qw+d_w\qw)}}
\\&
=
\E X\vw+\gl\vw(e^{-\gl\vw}-1)
  +O\lrpar{\gl\vw^2N\qw+\gl\vw^3e^{\gl\vw}(d_v\qw+d_w\qw)}
\\
&=
\E X\vw+\gl\vw e^{-\gl\vw}-\gl\vw
 +O\lrpar{(d_v^2d_w^3+d_v^3d_w^2)N^{-3}}
.  \end{split}
\end{equation*}
Summing these two equations we find
\begin{equation*}
  \begin{split}
\P(X\vw\le1)
=
(1+\gl\vw) e^{-\gl\vw}
 +O\lrpar{(d_v^2d_w^3+d_v^3d_w^2)N^{-3}}
  \end{split}
\end{equation*}
and the result follows, noting that $(1+\gl\vw) e^{-\gl\vw}$ is
bounded below since $\gl\vw=O(1)$.
\end{proof}

\section{Joint probabilities}\label{SC}

Our goal is to show that the indicators $I_u(\gx)$ and $J_e(\gx)$ are
almost independent for different $u$ and $e$; this is made precise in
the following lemma.

We define for convenience, for $u\in V_n$ and
$e=vw\in E_n$,
\begin{equation}\label{c0}
 \mu_u\=d_u^2/N 
\qquad \text{and} \qquad
\mu_e\=d_vd_w/N.
\end{equation}
It follows easily from \eqref{magnus} and a similar calculation for loops
that
\begin{equation}\label{cx}
\E \bigpar{X_u(\gx)\fall k}\le \mu_u^k
\qquad \text{and} \qquad
\E \bigpar{X_e(\gx)\fall k} \le \mu_e^k,
\qquad k\ge1.
\end{equation}
In particular, omitting the argument $\gx$,
\begin{equation}\label{cy}
\begin{aligned}
\P(I_u=1)&=\E I_u\le\E X_u\le\mu_u, 
\\
\P(J_e=1)&=\E J_e \le \E X_e\fall2\le\mu_e^2.
\end{aligned}  
\end{equation}
More precisely,
it follows easily from \refL{L1} that (for $N$ large at least)
$\P(I_u=1)=\Theta(\mu_u)$ and $\P(J_e=1)=\Theta(\mu_e^2)$ provided
$d_v,d_w\ge2$; this may help understanding our
estimates but will not be used below.

\begin{lemma}
  \label{L2}
Suppose \dbound. Let $l\ge0$ and $m\ge0$ be fixed.
For any
sequences of distinct vertices $u_1,\dots,u_l\in V_n$ and edges
$e_1,\dots,e_m\in E_n$, let $e_j=v_jw_j$ and let $F$ be the set of
vertices that appear at least twice in the list
$u_1,\dots,u_l,v_1,w_1,\dots,v_m,w_m$. Then, 
\begin{multline}\label{c1}
\E\biggpar{\prodil I_{u_i}(\gx) \prodjm J_{e_j}(\gx)}
=
\prodil\E( I_{u_i}(\gx)) \prodjm \E(J_{e_j}(\gx))
\\
+O\biggpar{\Bigpar{N\qw+\sum_{v\in F}d_v\qw}
\prodil\mu_{u_i} \prodjm \mu_{e_j}^2}
.  
\end{multline}
\end{lemma}
The implicit constant in the error term may depend on $l$ and $m$ but
not on $(u_i)_i$ and $(e_j)_j$. All similar statements below are to be
interpreted similarly.

The proof of \refL{L2} is long, and contains several other lemmas.
The idea of the proof 
is to use induction in $l+m$. In
the inductive step we select one of the indicators, $J_{e_1}$ say, and
then show that the product of the other indicators is almost
independent of $X_{e_1}$, and thus of $J_{e_1}$. In order to do so, we
would like to condition on the value of $X_{e_1}$. But the effects of
conditioning on $X_{e_1}=k$ are complicated and we find it difficult
to argue directly with these conditionings (see \refR{Rbad}).
Therefore, we begin with another, related but technically much simpler
conditioning. 

Fix two distinct vertices $v$ and $w$.
For $0\le k\le \min(d_v,d_w)$, let $\cE_k$ be
the event that the random configuration contains the $k$ pairs 
of half-edges
\set{\hv{i},\hw{i}}, $1\le i\le k$, 
and let the corresponding random
multigraph, \ie, $\gx$ conditioned on $\cE_k$, be denoted $\gx_k$.
$\gx_k$ thus contains at least $k$ edges between $v$ and $w$,
but there may be more.
Note that $\gx_0=\gx$.

We begin with an estimate related to \refL{L2}, but cruder.
\begin{lemma}
  \label{LD}
Suppose \dbound. 
Let $l$, $m$ and $r_1,\dots,r_l,s_1,\dots,s_m$ be fixed non-negative integers.
For any
sequences of distinct vertices $u_1,\dots,u_l\in V_n$ and edges
$e_1,\dots,e_m\in E_n$, 
\begin{equation}\label{d1}
\E\biggpar{\prodil X_{u_i}(\gx)\fall{r_i} \prodjm X_{e_j}(\gx)\fall{s_j}}
=
O\biggpar{\prodil\mu_{u_i}^{r_i} \prodjm \mu_{e_j}^{s_j}}
.  
\end{equation}
In particular,
\begin{equation}\label{d2}
\E\biggpar{\prodil I_{u_i}(\gx) \prodjm J_{e_j}(\gx)}
=
O\biggpar{\prodil\mu_{u_i} \prodjm \mu_{e_j}^{2}}
.  
\end{equation}

The estimates \eqref{d1} and \eqref{d2} hold with $\gx$ replaced by
$\gxk$ too, uniformly in $k$, provided the edges $e_1,\dots,e_m$ are
distinct from the edge $vw$ used to define $\gxk$. If $vw$ equals
some $e_j$, then \eqref{d1} still holds for $\gxk$, if we replace 
$X_{e_j}$  by $X_{e_j}-k$ when $e_j=vw$.
\end{lemma}
\begin{proof}
  We argue as for \eqref{magnus}. Let, again, $e_j=v_jw_j$ and 
let $t:=r_1+\dots+r_l+s_1+\dots+s_m$.
The expectation in \eqref{d1} is the number of $t$-tuples of disjoint
pairs of half-edges such that the first $r_1$ pairs have both
half-edges belonging to $u_1$, and so on, until the last $s_m$ that
each consist of one half-edge at $v_m$ and one at $w_m$, 
times the probability that a given such $t$-tuple is
contained in a random configuration. The number of such $t$-tuples is
at most $\prodil d_{u_i}^{2r_i}\prodjm (d_{v_j}d_{w_j})^{s_j}$
and the probability is $((2N-1)\dotsm(2N-2t+1))\qw<N^{-t}$
(provided $N\ge 2t$). The estimate \eqref{d1} follows, recalling
\eqref{c0}, and \eqref{d2} is
an immediate consequence since $I_u\le X_u$ and $J_e\le X_e\fall2$.

The same argument proves the estimates for $\gxk$. There is a minor
change in the probability above, replacing $N$ by $N-2k$; 
nevertheless, the estimates are uniform in $k$ because 
$k\le d_v=O(N\qq)$. (There may also be some $t$-tuples that are
excluded because they clash with the special pairs
\set{\hv{i},\hw{i}}, $i=1,\dots,k$; this only helps.)
\end{proof}

Let
$u_1,\dots,u_l\in V_n$ and $e_1,\dots,e_m\in E_n$
be as in \refL{L2}, and assume that $m\ge1$.
We choose $v=v_1$ and $w=w_1$, so $e_1=vw$, for the definition of $\gxk$.

If $k\ge1$, we can couple $\gxk$ and $\gxki$ as follows.
Start with a random configuration containing the $k$ special pairs 
\set{\hv{i},\hw{i}}. Then select, at random, a half-edge $x$ among all
half-edges except $\hv1,\dots,\hv{k},\hw1,\dots,\hw{k-1}$.
If $x=\hw{k}$ do nothing. Otherwise, let $y$ be the half-edge paired
to $x$; remove the two pairs \set{\hv{k},\hw{k}} and \set{x,y}
from the configuration
and replace them by \set{\hv{k},x} and \set{\hw{k},y}.
(This is called a switching; see McKay \cite{McKay} and McKay and
Wormald \cite{McKayWmoderate,McKayWo} for different but related arguments
with switchings.)

It is clear that this gives a configuration in $\cE_{k-1}$ with the
correct uniform distribution.
Passing to the multigraphs, we thus obtain a coupling of $\gxk$ and
$\gxki$ such that the two multigraphs differ (if at all) in that one
edge between $v$ and $w$ and one other 
edge have been deleted, and two new edges are added, one at $v$ and
one at $w$.

Let $Z$ denote the product 
$\prodil I_{u_i} \prodjmx2 J_{e_j}$ of the chosen indicators
except $J_{e_1}$. Define $F_1\subseteq\set{v,w}$ to be the set of
endpoints of $vw=e_1$ that also appear as some $u_i$ or as an
end-point of some other $e_j$; thus $F_1=F\cap\set{v,w}$. We claim the
following. 

\begin{lemma}
  \label{L3}
Suppose \dbound. 
With notations as above, 
uniformly in $k$ with $1\le k\le \min(d_v,d_w)$,
\begin{equation}\label{l3}
\E \bigpar{Z(\gxk)}-\E \bigpar{Z(\gxki)}
=
O\biggpar{\Bigpar{N\qw+\sum_{v\in F_1}d_v\qw}
\prodil\mu_{u_i} \prodjmx2 \mu_{e_j}^2}
.  
\end{equation}
\end{lemma}
\begin{proof}
We use the coupling above. Recall that $Z=0$ or $1$, so if $Z(\gxk)$
and $Z(\gxki)$ differ, then one of them equals 0 and the other equals
1. 

First, if $Z(\gxk)=1$ and $Z(\gxki)=0$, then the edge $xy$ deleted
from $\gxk$ must be either the only loop at some $u_i$, or one of
exactly two edges between $v_j$ and $w_j$ for some $j\ge2$.
Hence, for any configuration with $Z(\gxk)=1$, there are less than
$l+2m$ such edges, and the probability that one of them is destroyed
is less than $(l+2m)/(N-k) = O(1/N)$.
Hence,
\begin{equation}\label{c7}
  \P\bigpar{Z(\gxk)>Z(\gxki)}=O\bigpar{\E(Z(\gxk))/N}.
\end{equation}
Define $M\=\prodil\mu_{u_i} \prodjmx2 \mu_{e_j}^2$. 
By \refL{LD}, $\E Z(\gxk)=O(M)$, so the probability in \eqref{c7} is
$O(M/N)$, which is dominated by the \rhs{} of \eqref{l3}.

In the opposite direction, $Z(\gxk)=0$ and $Z(\gxki)=1$ may happen in
several ways. We list the possibilities as follows. 
(It is necessary but not necessarily sufficient for $Z(\gxk)<Z(\gxki)$
that one of them holds.)
\begin{romenumerate}
  \item
$v$ is an endpoint of one of the edges
$e_2,\dots,e_m$, say $v=v_2$ so $e_2=vw_2$;
the new edge from $v$ goes to $w_2$;
there already is (exactly) one edge between $v$ and
$w_2$ in $\gxk$; 
if we write $Z'=\prodil I_{u_i} \prodjmx3 J_{e_j}$, 
so that $Z=J_{e_2}Z'$, then $Z'(\gxk)=1$. 
\item
$v$ equals one of $u_1,\dots,u_l$, say $v=u_1$;
the new edge from $v$ is a loop;
if we write $Z'=\prodilx2 I_{u_i} \prodjmx2 J_{e_j}$, so that
$Z=I_{u_1}Z'$, 
then $Z'(\gxk)=1$. 
\item Two similar cases with $v$ replaced by $w$.
\item Both $v$ and $w$ are endpoints of edges $e_j$, say $v=v_2$ and
$w=w_3$, so that $e_2=vw_2$ and $e_3=wv_3$;
the two new edges go from $v$ to $w_2$ and from $w$ to $v_3$;
there are already such edges in $\gxk$;
if $Z''=\prodil I_{u_i} \prodjmx4 J_{e_j}$, so that
$Z=J_{e_2}J_{e_3}Z''$, 
then $Z''(\tgxk)=1$, where $\tgxk$ is $\gxk$ with one edge
between $w_2$ and $v_3$ deleted.
\item Both $v$ and $w$ equal some $u_i$, say $v=u_1$ and $w=u_2$;
the new edges from $v$ and $w$ are loops; 
if $Z''=\prodilx3 I_{u_i} \prodjmx2 J_{e_j}$, so that
$Z=I_{u_1}I_{u_2}Z''$, 
then $Z''(\gxk)=1$.
\item A similar mixed case where, say $v=v_2$ and $w=u_1$.
\item The same with $v$ and $w$ interchanged.
\end{romenumerate}

Consider case (i). For any configuration, the probability that the new
edge from $v$ goes to $w_2$ is $d_{w_2}/(2N-2k+1)=O(d_{w_2}/N)$. Since
we also need $Z'(\gxk)=1$ and $X_{e_2}(\gxk)\ge1$, the probability of case
(i) is at most
\begin{equation*}
  O\bigpar{d_{w_2}N\qw\E\bigpar{X_{e_2}(\gxk) Z'(\gxk)}}.
\end{equation*}
Now, by \refL{LD}, for convenience omitting the arguments $\gxk$ here
and often below in this proof,
\begin{equation*}
\E\bigpar{X_{e_2} Z'}
\le
\E\biggpar{\prodil X_{u_i}X_{e_2} \prodjmx3 X_{e_j}\fall2}
=O\biggpar{\prodil \mu_{u_i}\mu_{e_2} \prodjmx3\mu_{e_j}^2}
=O\bigpar{M/\mu_{e_2}}.
\end{equation*}
Moreover, $d_{w_2}/N=\mu_{e_2}/d_v$, so the probability of case (i) is
$O(M/d_v)$; note that the case only can happen if $v\in F_1$, so this
is covered by the \rhs{} of \eqref{l3}.

Case (ii) is similar (and slightly simpler).

Case (iii) occurs, by symmetry, with probability $O(M/d_w)$, and only
if $w\in F_1$.

In case (iv), the other destroyed edge must go between $w_2$ and
$v_3$. For any configuration, the probability that such an edge is
chosen is $O(X_{w_2v_3}/N)$. 
We study two subcases.
If one of the edges $e_4,\dots,e_m$
equals $w_2v_3$, say $e_4=w_2v_3$, then we, moreover, need at least
three edges between $w_2$ and $v_3$ in $\gxk$, since one of them is
destroyed. We also need $X_{e_2}(\gxk)\ge1$ and $X_{e_3}(\gxk)\ge1$.
Thus the probability of this case then is
\begin{equation*}
O\Bigpar{N\qw \E\bigpar{X_{w_2v_3}X_{e_2}X_{e_3} Z''(\tgxk)}}
= 
O\Bigpar{N\qw \E\bigpar{X_{e_2}X_{e_3}X_{e_4}\ett{X_{e_4}\ge3}	Z''}}.
\end{equation*}
By \refL{LD} we have
\begin{equation*}
  \begin{split}
\E\bigpar{X_{e_2}X_{e_3}X_{e_4}\ett{X_{e_4}\ge3}	Z''}
&\le
\E\biggpar{\prodil X_{u_i}\cdot X_{e_2}X_{e_3}X_{e_4}\fall3 
  \prodjmx5 X_{e_j}\fall2}
\\&
=O\biggpar{\prodil \mu_{u_i}\cdot\mu_{e_2}\mu_{e_3}\mu_{e_4}^3 
  \prodjmx5\mu_{e_j}^2}
\\&
=O\bigpar{M\mu_{e_4}/(\mu_{e_2}\mu_{e_3})}.
  \end{split}
\end{equation*}
In this case we have $\mu_{e_2}\mu_{e_3}=\mu_{e_1}\mu_{e_4}$,
so the probability is $O(N\qw M/\mu_{e_1})=O(M/d_vd_w)$.

In the second subcase, $w_2v_3$ does not equal any of $e_4,\dots,e_m$.
We then obtain similarly the probability
\begin{equation*}
  \begin{split}
O\Bigpar{N\qw\E\bigpar{X_{w_2v_3}X_{e_2}X_{e_3}Z''}}
&
=O\biggpar{N\qw\prodil \mu_{u_i}\cdot\mu_{w_2v_3}\mu_{e_2}\mu_{e_3}
  \prodjmx4\mu_{e_j}^2}
\\&
=O\bigpar{N\qw M\mu_{w_2v_3}/(\mu_{e_2}\mu_{e_3})},
  \end{split}
\end{equation*}
which again equals $O(M/(d_vd_w))$. 
Finally, note that in case (iv), $F_1=\set{v,w}$.

In case (v), the other destroyed edge is also an edge between $v$ and
$w$; given a configuration, the probability of this is 
$O((X_{vw}-k)/N)$. 
The probability of case (v) is thus
\begin{equation*}
  \begin{split}
O\Bigpar{N\qw\E\bigpar{(X_{vw}-k)Z''}}
&
=O\biggpar{N\qw\mu_{vw}\prodilx3 \mu_{u_i}  \prodjm\mu_{e_j}^2}
\\&
=O\bigpar{ M/(d_vd_w)}.
  \end{split}
\end{equation*}
$F_1=\set{v,w}$ in case (v) too.

Cases (vi) and (vii) are similar to case (iv), and lead to the same
estimate. Again $F_1=\set{v,w}$.

By \eqref{c7} and our estimates for the different cases above, the
probability that $Z(\gxk)$ and $Z(\gxki)$ differ is bounded by the
\rhs{} of \eqref{l3}, which completes the proof.
\end{proof}

We can now estimate the expectation of $Z(\gxk)$ conditioned on the
value of $X_{vw}(\gxk)$. We state only the result we need.
(See also \eqref{david}. These results can be rewritten as estimates of
conditional expectations.)

\begin{lemma}
  \label{L4}
Suppose \dbound. 
With notations as above, 
\begin{multline}\label{l4}
\E \bigpar{Z(\gx)J_{e_1}(\gx)}
\\
=
\E \bigpar{Z(\gx)} \E \bigpar{J_{e_1}(\gx)}
+
O\biggpar{\Bigpar{N\qw+\sum_{v\in F_1}d_v\qw}
\prodil\mu_{u_i} \prodjmx1 \mu_{e_j}^2}
.  
\end{multline}
\end{lemma}

\begin{proof}
  We can write $X_{vw}(\gx)\fall k=\sum_{\ga\in\cA} I_\ga$, where
  $\cA$ is the set of all ordered $k$-tuples of disjoint pairs $(x,y)$
  of half-edges with $x$ belonging to $v$ and $y$ to $w$, and $I_\ga$
  is the indicator that the $k$ pairs in $\ga$ all belong to the
  configuration.
By symmetry, $\E\bigpar{Z(\gx)\mid I_\ga=1}$ is the same for all
  $\ga\in\cA$; since $\gxk$ is obtained by conditioning $\gx$ on
  $I_\ga$ for a specific $\ga$, we thus have
$\E\bigpar{Z(\gx)\mid I_\ga=1} = \E\bigpar{Z(\gxk)}$ for all $\ga$.
Consequently,
\begin{equation}\label{erika}
  \begin{split}
\E\lrpar{X_{vw}(\gx)\fall k Z(\gx)}	
&=
\sum_{\ga\in\cA} \E\bigpar{I_\ga Z(\gx)}
=
\sum_{\ga\in\cA} \E\bigpar{ Z(\gx)\mid I_\ga=1}\P(I_\ga=1)
\\&
=
\E\bigpar{Z(\gxk)} \sum_{\ga\in\cA} \E I_\ga
=
\E\bigpar{Z(\gxk)}  \E \bigpar{X\vw(\gx)\fall k}.
  \end{split}
\end{equation}

We write the error term on the \rhs{} of \eqref{l3} as $O(R)$.
Since \eqref{l3} is uniform in $k$, and $\gx_0=\gx$, \refL{L3} yields
\begin{equation}\label{per}
\E\bigpar{Z(\gxk)}  
= \E\bigpar{Z(\gx)}  + O(kR).
\end{equation}
We now use \refL{LP}(ii) and (i) and find, for any $j$, using
\eqref{erika} and \eqref{per},
\begin{equation*}
  \begin{split}
  \E\bigpar{Z(\gx&)\cdot\ett{X\vw(\gx)=j}} 
= \sum_{k=j}^\infty (-1)^{k-j}\binom kj \frac{1}{k!}
  \E\bigpar{X\vw(\gx)\fall k Z(\gx)}
\\&
= \sum_{k=j}^\infty (-1)^{k-j}\binom kj \frac{1}{k!}
  \E\bigpar{X\vw(\gx)\fall k}\bigpar{\E Z(\gx)+O(kR)}
\\&
= \E\bigpar{Z(\gx)}\P\bigpar{X\vw(\gx)=j}
+O\lrpar{\sum_{k=j}^\infty \binom kj 
\frac{1}{k!}\E\bigpar{X\vw(\gx)\fall k}kR} .
 \end{split}
\end{equation*}
By \eqref{cx},
the sum inside the last $O$ is at most
\begin{equation*}
\sum_{k=j}^\infty \binom kj 
\frac{k}{k!}\mu\vw^k R
=  \sum_{l=0}^\infty \frac{j+l}{j!\,l!} \mu\vw^{j+l}R
= \biggpar{\frac{\mu\vw^j}{(j-1)!}+\frac{\mu\vw^{j+1}}{j!}}e^{\mu\vw}R.
\end{equation*}
Since $\mu\vw=O(1)$ by \eqref{dmaxbound}, we thus find, uniformly in $j\ge1$,
\begin{multline}
%\begin{equation} 
\label{david}
\E \bigpar{Z(\gx)\cdot\ett{X_{vw}(\gx)=j}}
\\
=
\E \bigpar{Z(\gx)}\P\bigpar{ X_{vw}(\gx)=j}
+
O\bigpar{\mu\vw^jR/(j-1)!}
, 
%\end{equation}
\end{multline}
which by summing over $j\ge2$ yields, again using $\mu\vw=O(1)$ and
recalling that $vw=e_1$,
\begin{equation*}
\E \bigpar{Z(\gx)J_{e_1}(\gx)}
=
\E \bigpar{Z(\gx)} \E\bigpar{J_{e_1}(\gx)}
+
O\bigpar{\mu_{e_1}^2R}
,
\end{equation*}
the sought estimate.
\end{proof}

\begin{proof}[Proof of \refL{L2}] 
As said above, we use induction on $l+m$. The result is trivial if
$l+m=0$ or $l+m=1$. If $m\ge1$, we use \refL{L4}; the result follows
from \eqref{l4} together with the induction hypothesis applied to
$\E(Z(\gx))$ and the estimate $\E\bigpar{J_{e_1}(\gx)}=O(\mu_{e_1}^2)$
from \refL{LD}.

If $m=0$, we study a product $\prodil I_{u_i}$ of loop indicators
only. We then modify the proof above, using loops instead of multiple
edges in the conditionings. More precisely,
we now let $\gxk$ be $\gx$ conditioned on the configuration
containing the $k$ specific pairs $(\hu{2i-1},\hu{2i})$,
$i=1,\dots,k$, of half-edges at $u$. We couple $\gxk$ and $\gxki$ as
above (with obvious modifications). 
In this case, the switching from $\gxk$ to $\gxki$ cannot create any new
loops. Hence, if $Z\=\prodilx2 I_{u_i}$, we have $Z(\gxk)\ge
Z(\gxki)$.
We obtain \eqref{c7} exactly as before, and since \refL{LD} still
holds, this shows that \refL{L3} holds, now with $F_1=\emptyset$ and the
error term $O(N\qw\prodilx2\mu_{u_i})$.
It follows that \refL{L4} holds too (with $J_{e_1}$ replaced by
$I_{u_1}$ and $F_1=\emptyset$) by the same proof as above.
This enables us to complete the induction step in the case $m=0$ too.
\end{proof}

\begin{remark}
  Similar arguments show that Lemmas \refand{L3}{L4}, with obvious
  modifications, hold in this setting, where we condition on loops at
  $u$, also for $m>0$. A variation of our proof of \refL{L2} would be
  to use this as long as $l>0$; the result in our \refL{L3} then is
  needed only when $l=0$, which eliminates cases (ii), (v), (vi),
  (vii) from the proof. On the other hand, we have to consider new
  cases for the loop version, so the total amount of work is about the same.
\end{remark}

\begin{remark}
  When conditioning on loops, it is possible to argue directly with
  conditioning on $X_u=k$, using a coupling similar to the one for
  $\gxk$ above; we thus do not need the detour with $\gxk$ and
  \refL{LP} used above. However, as said above, 
in order to treat multiple edges, the method used here
seems to be much simpler. A possible alternative would be to use the
methods in McKay \cite{McKay} and 
McKay and Wormald \cite{McKayWmoderate,McKayWo};
we can interpret the arguments there as showing that suitable
switchings yield an approximate, but not exact, coupling
when we condition on exact numbers of edges in different positions.
\end{remark}

\begin{remark}
  \label{Rbad}
A small example that illustrates some of the complications when
conditioning on a given number of edges between two vertices is
obtained by taking three vertices $1,2,3$ of degree 2 each.
Note that if $X_{12}=1$, then the multigraph must be a cycle; in
particular, $X_3=0$. On the other hand, $X_3=1$ is possible for
$X_{12}=0$; this shows that it is impossible to couple the multigraphs
conditioned on $X_{12}=0$ and on $X_{12}=1$ 
by moving only two edges as in the proof above.
Note also that $X_3=0$ is possible also when
$X_{12}=2$; there is thus a surprising non-convexity.
\end{remark}

\section{The proofs are completed}\label{Sproofs}

\begin{proof}[Proof of \refT{T2}]
We begin by observing that the two expressions given in \eqref{t2a}
and \eqref{t2b} are
equivalent. Indeed, if we define $\gL$ by \eqref{Lambda}, then
\begin{equation*}
  \begin{split}
\hskip 3em&\hskip-3em
\tfrac12\sumi\gl\ii +\tfrac12\sum_{i<j}\gl\ij^2
=
\tfrac12\sumi\gl\ii 
 +\tfrac14\Bigpar{\sum_{i,j}\gl\ij^2 - \sum_{i}\gl\ii^2}
\\&
=	
\frac12\sumi\frac{d_i(d_i-1)}{2N} 
+ \frac14\Bigpar{\sum_{i,j}\frac{d_i(d_i-1)d_j(d_j-1)}{4N^2}
 -\sumi\frac{d_i^2(d_i-1)^2}{4N^2}}
\\&
=\gL+\gL^2-\frac{\sumi d_i^2(d_i-1)^2}{16 N^2}
  \end{split}
\end{equation*}
and
\begin{equation}\label{julie}
  \gL^2+\gL 
= \Bigpar{\gL+\frac12}^2-\frac14
= \frac14\Bigpar{\frac{\sumi d_i^2}{2 N}}^2 -\frac14.
\end{equation}

We note for future reference that
\begin{equation*}
\sumi\gl\ii 
=
\frac1{2N}\sumi(d_i^2-d_i)=\tfrac12N\qw\sumdd-1
\end{equation*}
and 
$\gl-\log(1+\gl)\ge0$ when
$\gl\ge0$, and thus the \rhs{} of \eqref{t2a} can be estimated from above by
\begin{equation}  \label{j1}
\text{\rhs{} of (\ref{t2a})}
\le \exp\biggpar{\frac12-\frac{\sumdd}{4N}}+o(1).
\end{equation}
Similarly,
$\log(1+\gl)-\gl+\tfrac12\gl^2\ge0$ when
$\gl\ge0$, and thus 
\begin{equation}  \label{j2}
\text{\rhs{} of (\ref{t2b})}
\ge \exp\biggpar{-\Bigpar{\frac{\sumdd}{4N}}^2} +o(1).
\end{equation}
In particular, since we just have shown that these two \rhs{s} are the
same, they tend to 0 if and only if $\sumdd/N\to\infty$.

Next, suppose first that $\sumdd=O(N)$.
Recall $Y(\gx)=\sum_{u\in V_n}I_u(\gx)+\sum_{e\in E_n}J_e(\gx)$
defined in \eqref{y}.
As said above, \refL{L2} shows that the random variables $I_u(\gx)$
and $J_e(\gx)$ are almost independent. We can compare them with truly
independent variables as follows.

Let $\bIu$ and $\bJe$ be \emph{independent} 0--1 valued
random variables such that 
$\P(\bIu=1)=\P\bigpar{I_u(\gx)=1}$
and 
$\P(\bJe=1)=\P\bigpar{J_e(\gx)=1}$,
and let
\begin{equation*}
\bY\=\sum_{u\in V_n}\bI_u+\sum_{e\in E_n}\bJ_e.
\end{equation*}

Fix $k\ge1$. We use \refL{L2} for all pairs $(l,m)$ with $l+m=k$ and
sum \eqref{c1} over all such $(l,m)$ and all distinct $u_1,\dots,u_l$
and $e_1,\dots,e_m$, 
multiplying by the symmetry factor $\binom kl$, and
noting that the first term on the \rhs{} of
\eqref{c1} can be written $\E\bigpar{\prod_i\bI_{u_i}\prod_j\bJ_{e_j}}$.
This gives
\begin{equation}
  \label{e1}
\E\bigpar{Y(\gx)\fall k}
=
\E(\bY\fall k)
+O\biggpar{\sum\biggpar{\Bigpar{N\qw+\sum_{v\in F}d_v\qw}
\prodil\mu_{u_i} \prodjm \mu_{e_j}^2}},
\end{equation}
summing over all such $l,m,(u_i)_i,(e_j)_j$ and with $F$ depending on
them as in \refL{L2}.

Consider one term in the sum in \eqref{e1},
write as usual $e_j=v_jw_j$, 
and let $H$ be the multigraph with vertices $V(H)=\set{u_i}\cup\set{v_j,w_j}$
 and one loop at each $u_i$ and two parallel edges between $v_j$ and
 $w_j$ for each $j\le m$. Let $\dvh$ be the degree of vertex $v$ in
 $H$ and note that each degree $\dvh$ is even, and thus at least 2,
 and that $F$ is the set of vertices with $\dvh\ge4$. We have
 \begin{equation}\label{e2}
\prodil\mu_{u_i} \prodjm \mu_{e_j}^2 
=N^{-e(H)}\prod_{v\in V(H)} d_v^{\dvh},   
 \end{equation}
where $e(H)=l+2m$ is the number of edges in $H$.

We group the terms in the sum in \eqref{e1} according to the
isomorphism type of $H$. 
Fix one such type $\cH$, and let it have $h$ vertices with degrees
$a_1,\dots,a_h$ (in some order) and $b$ edges; thus
$b=\tfrac12\sum_{j=1}^h a_j$.
The corresponding $H$ are obtained by selecting vertices
$v_1,\dots,v_h\in V_n$; these have to be distinct and it may happen
that some permutations give the same $H$, but we ignore this, thus
overcounting, and obtain from \eqref{e2}
that 
  \begin{multline}\label{e2x}
\sum_{H\in\cH}\biggpar{\prodil\mu_{u_i} \prodjm \mu_{e_j}^2}
\le
  \sum_{v_1,\dots,v_h\in V_n} \biggpar{N^{-b} \prod_{j=1}^h d_{v_j}^{a_j}}
=
 N^{-b} \prod_{j=1}^h  \biggpar{\sum_{v_j\in V_n}d_{v_j}^{a_j}}
\\
=O\lrpar{N^{-b+\sum_j a_j/2}}
=O(1)
  \end{multline}
by \eqref{sofie}, since each $a_j\ge2$.

Furthermore, let $\hh:=\set{i\in\set{1,\dots,h}: a_i\ge 4}$. Thus, if $H$ is
obtained by choosing vertices $v_1,\dots,v_h\in V_n$, then
$F=\set{v_i:i\in\hh}$. Consequently,
  \begin{multline*}
\sum_{H\in\cH}\biggpar{{\sum_{v\in F}d_v\qw} \prodil\mu_{u_i} 
  \prodjm \mu_{e_j}^2}
\le
\sum_{v_1,\dots,v_h\in V_n}\sum_{i\in\hh}
  \biggpar{ d_{v_i}\qw N^{-b} \prod_{j=1}^h d_{v_j}^{a_j}}
\\
=
\sum_{i\in\hh}\biggpar{
 N^{-b} \prod_{j=1}^h  \sum_{v_j\in V_n}d_{v_j}^{a_j-\gd{ji}}}
=O\lrpar{N^{-b+\sum_j a_j/2-1/2}}
=O(N^{-1/2}),
  \end{multline*}
by \eqref{sofie}, since each $a_i\ge4$ if $i\in\hh$ and thus
$a_j-\gd_{ij}\ge 2$ for every $j$.

Combining this with \eqref{e2x}, we see that the sum in \eqref{e1}, 
summing only over $H\in\cH$, is $O(N^{-1/2})$.
There is only a finite number of isomorphism types $\cH$
for a given $k$, and thus we obtain the same estimate for the full
sum. Consequently, \eqref{e1} yields
\begin{equation}
\E\bigpar{Y(\gx)\fall k}
=
\E(\bY\fall k)+ O(N^{-1/2}),
\end{equation}
for every fixed $k$.

We use \refL{LB} with $\bY$ and $Y(\gx)$ (in this order). 
We have just verified \eqref{lb2}.
To verify
\eqref{lb1} we take $R=2$ 
(any $R<\infty$ would do by a similar argument)
and find, using \eqref{cy} and \eqref{dbound}
\begin{equation*}
  \begin{split}
\E\Bigpar{2^{\bY}}
&=  \prod_{u\in V_n}\E2^{\bIu}\prod_{e\in E_n}\E2^{\bJe}	
=  \prod_{u\in V_n}(1+\E\bIu)\prod_{e\in E_n}(1+\E\bJe)	
\\&
\le  \prod_{u\in V_n}\exp\xpar{\E\bIu}\prod_{e\in E_n}\exp\xpar{\E\bJe}	
= \exp\biggpar{\sum_{u\in V_n}\E\bIu+\sum_{e\in E_n}\E\bJe	}
\\&
\le\exp\biggpar{\sum_{u\in V_n}\frac{d_u^2}{N}
 +\sum_{vw\in E_n}\frac{d_v^2d_w^2}{N^2}}
=O(1).
  \end{split}
\end{equation*}
Consequently, \refL{LB} applies and yields
\begin{equation*}
  \begin{split}
\P\bigpar{\gx\simple}	
&=
\P(Y(\gx)=0)
=
\P(\bY=0)+o(1)
\\&
=
\prod_{u\in V_n}\P(\bIu=0) \prod_{e\in E_n}\P(\bJe=0) +o(1)
\\&
=
\exp\biggpar{\sum_{u\in V_n}\log\P(\bIu=0)+ \sum_{e\in E_n}\log\P(\bJe=0)} 
+o(1).
  \end{split}
\end{equation*}
Furthermore, \refL{L1} yields
\begin{multline*} 
%\begin{equation*}
%  \begin{split}
-\biggpar{\sum_{u\in V_n}\log\P(\bIu=0)+ \sum_{e\in
 E_n}\log\P(\bJe=0)} 
\\
=\sum_{u\in V_n}\frac{d_u(d_u-1)}{4N}
+\sum_{vw\in E_n}\bigpar{-\log(1+\gl\vw)+\gl\vw}
\\
+O\parfrac{\sum_{v}d_v^3}{N^2}
+O\parfrac{\sum_v d_v^3\sum_w d_w^2}{N^3},	
\end{multline*}
%  \end{split}
%\end{equation*}
where the two error terms ore $O(N\qqw)$ by \eqref{sofie}.

This verifies \eqref{t2a} and thus \refT{T2} in the case $\sumdd=O(N)$.

Next, suppose that $\sumdd\to\infty$. Fix a number $A>2$. For all
large $N$ (or $\nu$ in the formulation of \refT{T1}), $\sumdd>AN$, so
we may assume this inequality.

Let $j$ be an index with $d_j>1$. We then may modify the sequence
$\dn$ by decreasing $d_j$ to $d_j-1$ and adding a new element
$d_{n+1}=1$. This means that we split one of the vertices in $\gx$
into two. Note that this splitting increases the number $n$ of
vertices, but preserves the number $N$ of edges. We can repeat and
continue splitting vertices (in arbitrary order) until all degrees
$d_i\le1$; then $\sumdd=\sumd=2N$.

Let us stop this splitting the first time $\sumdd\le AN$ and denote
the resulting sequence by \hdn. Thus $\sumhdd\le AN$.
Since we have assumed $\sumdd>AN$, we have performed at least one
split. If the last split was at $j$, the sequence preceding $\hdn$ is
$(\hd_i+\gd_{ij})_1^{\hn-1}$, and thus
\begin{equation*}
  AN < \sum_{i=1}^{\hn-1}(\hd_i+\gd_{ij})^2
=\sum_{i=1}^{\hn} \hd_i^2 +2\hd_j+1-1
\le\sum_{i=1}^{\hn} \hd_i^2 +2\sqrt{AN},
\end{equation*}
because $\hd_j^2\le\sumhdd\le AN$.
Consequently, 
\begin{equation*}
  AN-2\sqrt{AN} \le \sumhdd \le AN
\end{equation*}
and thus, in the limit, $\sumhdd/N\to A$.

Let $\hgx=\gx\bigpar{\hn,\hdn}$. Since $\sumhdd=O(N)$, we can by the
already proven part apply \eqref{t2a} to $\hgx$ and obtain, using \eqref{j1},
\begin{equation*}
\P(\hgx\simple)
\le \exp\biggpar{\frac12-\frac{\sumhdd}{4N}}+o(1)
= \exp\biggpar{\frac12-\frac{A}{4}}+o(1).
\end{equation*}
Furthermore, since $\hgx$ is constructed from $\gx$ by splitting
vertices, $\hgx$ is simple whenever $\gx$ is, and thus
$\P(\gx\simple)\le\P(\hgx\simple)$.
Consequently,
\begin{equation*}
  \limsup \P(\gx\simple)
\le \limsup\P(\hgx\simple)
\le\exp\bigpar{-(A-2)/4}.
\end{equation*}
Since $A$ can be chosen arbitrarily large, this shows that if
$\sumdd/N\to\infty$, then $\P(\gx\simple)\to0$.

Combined with \eqref{j1}, this shows that \eqref{t2a} holds when
$\sumdd/N\to\infty$ (with both sides tending to 0).

Finally, for an arbitrary sequence of sequences $\dn$, we can for
every subsequence find a subsequence where either $\sumdd/N=O(1)$
or $\sumdd/N\to\infty$, and thus \eqref{t2a} holds along the
subsequence by one of the two cases treated above. It follows 
by the subsubsequence principle that
\eqref{t2a} holds.
\end{proof}

\begin{proof}[Proof of \refT{T1}]
Part (ii) follows immediately from \refT{T2} and \eqref{j1},
\eqref{j2}.

If we apply this to subsequences, we see that
$\liminf\P(\gx_\nu\simple)>0$ if and only if there is no subsequence along
which $\sum_i(d_i\nnu)^2/N_\nu\to\infty$, which proves (i).
\end{proof}

\begin{proof}[Proof of \refC{C2}]
  The two expressions are equivalent  by \eqref{julie}.

If $\sumi d_i^2/N\to\infty$, then the \rhs{} tends to 0, and so does
the \lhs{} by \refT{T1}. Hence, by the subsubsequence principle,
it remains only to show the result assuming $\sumi d_i^2=O(N)$.
In this case we have
\begin{equation*}
  \frac{\sumi d_i^2(d_i-1)^2}{16 N^2}
\le \frac{(\maxd)^2}{N}\frac{\sumi d_i^2}{N}
=o(1)
\end{equation*}
and, since $\log(1+x)-x+\tfrac12x^2=O(x^3)$ for
$x\ge0$,
\begin{equation*}
  \sum_{i<j}\bigpar{\log(1+\gl\ij)-\gl\ij+\tfrac12\gl\ij^2}
=O\biggpar{\sum_{i,j}\gl\ij^3}
\end{equation*}
and
\begin{equation*}
\sum_{i,j}\gl\ij^3
\le
\sum_{i,j}\frac{d_i^3d_j^3}{N^{3}}
\le \frac{(\maxd)^2}{N}\biggpar{\frac{\sumi d_i^2}{N}}^2
=o(1).
\end{equation*}
Hence, in this case, the formula in \refC{C2} follows from \eqref{t2b}.
\end{proof}

\section{Poisson approximation}\label{Spoisson}

As remarked in the introduction, when $\maxd=o(N\qq)$, it is easy to
prove that \eqref{a3} implies \eqref{a2} by the Poisson approximation
method of Bollob\'as \cite{BB80,bollobas}.
Since this is related to the method above, but much simpler, and we
find it interesting to compare the two methods, we describe this
method here, thus obtaining an alternative proof of \refC{C2}.
We  assume $\sumdd=O(N)$ throughout this section.

The main idea is to study the random variable
\begin{equation}\label{ty}
\tY\=\sum_{u\in V_n}X_u+\sum_{e\in E_n}\binom{X_e}{2},
\end{equation}
which counts the number of loops and pairs of parallel edges 
(excluding loops)
in $\gx$
(we omit the argument $\gx$ in this section). Compare this with $Y$
defined in \eqref{y}, and note that
\begin{equation*}
  \gx\simple \iff \tY=0 \iff Y=0.
\end{equation*}

\begin{theorem}
  \label{Tpoisson}
Assume that $\Ntoo$, $\sumdd=O(N)$ and $\maxd=o(N\qq)$.
Let $\gL\=\frac1{2N}\sumin\binom{d_i}2$
as in \eqref{Lambda}.
Then
\begin{equation}\label{p2}
  \dtv\bigpar{\tY,\Po(\gL+\gL^2)}\to0,
\end{equation}
and thus
\begin{equation}\label{p2a}
  \begin{split}
\P\bigpar{\gndx\simple}
= \P(\tY=0)
=
\exp\bigpar{-\gL-\gL^2}+o(1).
  \end{split}
\end{equation}
\end{theorem}

If $\gL\to\gl$ for some $\gl\in[0,\infty)$, then \eqref{p2} is
equivalent to $\tY\dto\Po(\gl+\gl^2)$. (By the subsubsequence
principle, it suffices to consider this case.)

\begin{proof}[Sketch of proof]
We can write $\tY=\sum_{\ga\in\cA}I_\ga+\sum_{\gb\in\cB}J_\gb$,
where $\cA$ is the set of all pairs \set{\hu{i},\hu{j}} of
half-edges (correponding to loops),
and $\cB$ is the set of all pairs of pairs 
\set{\set{\hv{i},\hw{j}},\set{\hv{k},\hw{l}}} of distinct half-edges
(corresponding to pairs of parallel edges).

Thus, similarly to \eqref{magnus},
\begin{equation}\label{p3}
  \begin{split}
\E\tY
&=
\sum_{\ga\in\cA}\E I_\ga	 +\sum_{\gb\in\cB}\E J_\gb
\\&
=
\sum_{u\in V_n} \frac{d_u(d_u-1)}{2(2N-1)}
+
\frac12\sum_{v\neq w} \frac12\frac{d_v(d_v-1)d_w(d_w-1)}{(2N-1)(2N-3)}
\\&
=\gL+\gL^2+o(1).
  \end{split}
\end{equation}

Moreover, it is easy to compute the expectation of a product of the form
$\E\bigpar{I_{\ga_1}\dotsm I_{\ga_l}J_{\gb_1}\dotsm J_{\gb_m}}$;
it is just the probability that a random configuration contains all
pairs occurring in $\ga_1,\dots,\ga_l,\gb_1,\dots,\gb_m$.
If two of these pairs intersect in exactly one half-edge, the
probability is 0; otherwise it is $(2N)^{-b}(1+O(1/N))$, where $b$ is
the number of different pairs. (Note that we may have, for example,
$\gb_1=\set{\set{\hv{1},\hw{1}},\set{\hv{2},\hw{2}}}$ and
$\gb_2=\set{\set{\hv{1},\hw{1}},\set{\hv{3},\hw{3}}}$, with one pair in
common; thus $b\le l+2m$, but strict inequality is possible.)

We can compute factorial moments $\E\bigpar{\tY\fall k}$ by summing
such expectations of products with $l+m=k$.
For each term
$\E\bigpar{I_{\ga_1}\dotsm J_{\gb_m}}$, let $H$ be the multigraph, with
vertex set a subset of $V_n$, obtained by joining each pair occurring
in $\ga_1\dots,\gb_m$ (taking repeated pairs just once) to an edge, and then
deleting all unused (\ie, isolated) vertices in $V_n$.
It is easy to estimate the sum of these terms for a given $H$, and we
obtain
$O\bigpar{N^{-e(H)}\prod_{v\in V(H)} d_v^{\dvh}}$ as in \eqref{e2}.
As in the proof of \refT{T2}, we then group the terms according to the
isomorphism type $\cH$ of $H$. (There are more types $\cH$ now, but
that does not matter.)

Since now $\maxd=o(N\qq)$, \eqref{sofie} is improved to
\begin{equation}\label{sofiex}
  \sumv d_v^k = o(N^{k/2})
\end{equation}
for every fixed $k\ge3$, and it follows that the sum for a given $\cH$
is $o(1)$ as soon as $\cH$ has at least one vertex with degree
$\ge3$. The only remaining case is when $\cH$, and thus $H$, consists
of $l$ and $m$ vertex-disjoint loops and double edges; in this case
\begin{equation}
  \label{p4}
  \begin{split}
\E\Bigpar{\prodil I_{\ga_i}\prodjm J_{\gb_j}}
&=
\bigpar{(2N-1)\dotsm(2N-2l-4m+1)}\qw
\\&
=\bigpar{1+O(1/N)}\prodil\E I_{\ga_i}\prodjm\E J_{\gb_j}.	
  \end{split}
\end{equation}

Similarly, we can expand 
$(\E\tY)^k
=\bigpar{\sum_{\ga\in\cA}\E I_\ga+\sum_{\gb\in\cB}\E J_\gb}^k$
as a sum of terms $\prodil\E I_{\ga_i}\prodjm\E J_{\gb_j}$ with
$l+m=k$. (Now, repetitions are allowed among $\ga_i$ and $\gb_j$.)
If we introduce $H$ and $\cH$ as above, we see again that the sum of
all terms with a given $\cH$ is $o(1)$ except when $\cH$ consists of 
$l$ and $m$ vertex-disjoint loops and double edges. The terms
occurring in this case are the same as in \eqref{p4}, and hence their
sums differ by $O(1/N)$ only (since these sums are $O(1)$, see
\eqref{p3}).

Consequently, summing over all $\cH$ and using \eqref{p3}, for every
$k\ge1$,
\begin{equation*}
  \E \bigpar{\tY\fall k}
=  \bigpar{\E\tY}^k + o(1)=(\gL+\gL^2)^k+o(1).
\end{equation*}
If $\gL\to\gl$, this shows $\tY\dto \Po(\gl+\gl^2)$ by the method of
moments.
In general, we obtain \eqref{p2a} and \eqref{p2} by \refL{LB} and \refR{RB}.
\end{proof}

\begin{remark}
This argument further shows that, asymptotically, the
number of loops is $\Po(\gL)$ and the number of pairs of double edges
is $\Po(\gL^2)$, with these numbers asymptotically independent.
\end{remark}

In order to compare this method with the one in the preceding
sections, note that $\tY\ge Y$ and that $\tY=Y$ if and only if there
are no double loops or triple edges. It is easy to see
that if $\sumdd=O(N)$ and $\maxd=o(N\qq)$, 
then, using \eqref{cx} and \eqref{sofiex}, 
\begin{equation}\label{p5}
  \begin{split}
\P(\tY\neq Y) 
&\le \sum_{u\in V_n} \E (X_u\fall2) + \sum_{e\in E_n} \E(X_e\fall3)
\\
&= O\parfrac{\sum_u d_u^4}{N^2} + O\parfrac{\bigpar{\sum_v
	d_v^3}^2}{N^3}
=o(1),
  \end{split}
\end{equation}
so in this case the two variables are equivalent asymptotically. In
particular, \refT{Tpoisson} is valid for $Y$ too.
It is evident that the argument to estimate factorial moments of $\tY$
in the proof of \refT{Tpoisson} is
much shorter that the argument to estimate factorial moments of $Y$ in
the preceding sections. The reason for the difference is the ease with
which we can compute $\E(I_{\ga_1}\dotsm J_{\gb_m})$ for a random
configuration. Hence the proof of Theorem 7.1 is preferable in this case.

On the other hand, if $\max d_i=\Theta(N\qq)$, 
still assuming $\sumdd=O(N)$, 
there are several
complications. 
Let us for simplicity assume that $d_1\ge d_2\ge\dots$, and that
$d_1\sim c_1 N\qq$ with $c_1>0$. Then $X_1\dto\Po(c_1^2/4)$, so
$\lim\P(X_u>1)>0$ and \eqref{p5} fails. 

Moreover, \cf{} \eqref{p3},
\begin{equation*}
  \E \tY 
=\tfrac12\sum_i\gl\ii +\tfrac12\sum_{i<j}\gl\ij^2
+o(1),
\end{equation*}
so we can write \eqref{t2a} and \eqref{t2b} as
\begin{equation*}
%\begin{multline}
\P\bigpar{\gx\simple}
=
\exp\Bigpar{-\E\tY
+\sum_{i<j}\bigpar{\log(1+\gl\ij)-\gl\ij+\tfrac12\gl\ij^2}}
+o(1).
%\end{multline}  
\end{equation*}
Except in the case $d_2=o(N\qq)$, we cannot ignore the terms with
$\gl\ij$ in the exponent; if, say, $d_2\sim c_2 N\qq$ with $c_2>0$, then
$\gl_{12}\to c_1c_2/2>0$. Consequently, \refT{T2} shows that in this
case, $\P(\tY=0)=\P(\gx\simple)$ is not well approximated by
$\exp(-\E\tY)$,
which shows that $\tY$ is \emph{not} asymptotically Poisson distributed.
(The reason is terms like $\binom {X_{12}}2$ in \eqref{ty}, 
where $X_{12}\dto\Po(c_1c_2/2)$.)

Further, we have shown in \refS{Sproofs} that $Y$
asymptotically can be regarded as the sum $\bY$ of independent
indicators, 
but in this case  $\lim \P(I_1=1)>0$, and thus these indicators do not
all have small expectations; hence $Y$ is not asymptotically Poisson
distributed. 

Any attempt to show Poisson convergence of either $\tY$ or $Y$ is thus
doomed to fail unless $\maxd=o(N\qq)$. It seems difficult to find the
asymptotic distribution of $\tY$ directly; even if we could show that
the moments converge, the moments grow too rapidly for the method of
momemts to be applicable (at least with the Carleman criterion, see
Section 4.10 in \cite{Gut}). This is the reason for studying $Y$
above; as we have seen above, the distribution is asymptotically nice,
even if our proof is rather complicated.

\begin{ack}
I thank Bela Bollob\'as and Nick Wormald for helpful comments.  
\end{ack}

\newcommand\AAP{\emph{Adv. Appl. Probab.} }
\newcommand\JAP{\emph{J. Appl. Probab.} }
\newcommand\JAMS{\emph{J. \AMS} }
\newcommand\MAMS{\emph{Memoirs \AMS} }
\newcommand\PAMS{\emph{Proc. \AMS} }
\newcommand\TAMS{\emph{Trans. \AMS} }
\newcommand\AnnMS{\emph{Ann. Math. Statist.} }
\newcommand\AnnPr{\emph{Ann. Probab.} }
\newcommand\CPC{\emph{Combin. Probab. Comput.} }
\newcommand\JMAA{\emph{J. Math. Anal. Appl.} }
\newcommand\RSA{\emph{Random Structures Algorithms} }
\newcommand\ZW{\emph{Z. Wahrsch. Verw. Gebiete} }
\newcommand\DMTCS{\jour{Discr. Math. Theor. Comput. Sci.} }

\newcommand\AMS{Amer. Math. Soc.}
\newcommand\Springer{Springer-Verlag}
\newcommand\Wiley{Wiley}

\newcommand\vol{\textbf}
\newcommand\jour{\emph}
\newcommand\book{\emph}
\newcommand\inbook{\emph}
\def\no#1#2,{\unskip#2, no. #1,} %(typeset after year) 
\newcommand\toappear{\unskip, to appear}

\newcommand\webcite[1]{\hfil\penalty0\texttt{\def~{\~{}}#1}\hfill\hfill}
\newcommand\webcitesvante{\webcite{http://www.math.uu.se/\~{}svante/papers/}}
\newcommand\arxiv[1]{\webcite{http://arxiv.org/#1}}

\def\nobibitem#1\par{}


\begin{thebibliography}{99}

%\frenchspacings

\bibitem{BenderC}
E.~A. Bender \& E.~R. Canfield,
The asymptotic number of labeled graphs with given degree sequences.
\jour{J. Combin. Theory Ser. A}, \vol{24}\no3 (1978), 296--307.


\bibitem{BB80}
B. Bollob\'as, 
A probabilistic proof of an asymptotic formula for the number of
labelled regular graphs,  
\jour{European J. Comb.} \vol1 (1980), 311--316.

\bibitem{bollobas}
B. Bollob\'as, \book{Random Graphs}, 2nd ed., Cambridge Univ. Press,
Cambridge, 2001.

\bibitem{cooper-hyper}
C. Cooper,
The cores of random hypergraphs with a given degree sequence.  
\RSA \vol{25} \no4 (2004), 353--375.

\bibitem{cooperF-digraph}
C. Cooper \& A. Frieze, 
The size of the largest strongly connected component of a random
digraph with a given degree sequence.  
\CPC \vol{13} \no3  (2004), 319--337.

\bibitem{Gut}
A. Gut,
\book{Probability: A Graduate Course}.
Springer, New York, 2005.

\bibitem{SJ184}
S. Janson \& M. Luczak,
A simple solution to the $k$-core problem. 
\RSA\unskip, to appear.
\arxiv{math.CO/0508453}

\bibitem{JLR}
S. Janson, T. \L uczak \& A. Ruci\'nski,
\book{Random Graphs}.
\Wiley, New York, 2000.

\bibitem{McKayrect}
B. D. McKay, 
Asymptotics for $0$-$1$ matrices with prescribed line sums.  
\inbook{Enumeration and design (Waterloo, Ont., 1982)}, pp. 225--238, 
Academic Press, Toronto, ON, 1984.

\bibitem{McKay}
B. D. McKay, 
Asymptotics for symmetric $0$-$1$ matrices with prescribed row sums.  
\jour{Ars Combin.}  \vol{19A}  (1985),  15--25.

\bibitem{McKayWmoderate}
B. D. McKay \& N. C. Wormald, 
Uniform generation of random regular graphs of moderate degree.
\emph{J. Algorithms} \vol{11}\no1 (1990), 52--67.
05C80 (68R10)

\bibitem{McKayWhigh}
B. D. McKay \& N. C. Wormald, 
Asymptotic enumeration by degree sequence of graphs of high degree.  
\jour{European J. Combin.}  \vol{11}\no6  (1990),  565--580.

\bibitem{McKayWo}
B. D. McKay \& N. C. Wormald, 
Asymptotic enumeration by degree sequence of graphs with degrees
$o(n\sp {1/2})$.  
\jour{Combinatorica} \vol{11}\no4  (1991), 369--382.

\bibitem{WormaldPhD}
N. C. Wormald, 
Some problems in the enumeration of labelled graphs.
Ph. D. thesis, University of Newcastle, 1978.

\bibitem{Wormald81}
N. C. Wormald, 
The asymptotic distribution of short cycles in random regular graphs.  
\jour{J. Combin. Theory Ser. B}  \vol{31}\no2  (1981), 168--182.



\end{thebibliography}
\end{document}